\documentclass[12pt,a4paper]{amsart}
\usepackage{amsmath}
\usepackage{amssymb}

\usepackage{amsthm}
\usepackage{color}
\usepackage{enumerate}
\usepackage[dvips]{graphicx}
\usepackage{float}
\usepackage{wrapfig}
\usepackage{layout}
\usepackage{comment}

\newcommand{\T}{\mathbb{T}}

\DeclareMathOperator\Lip{Lip}

\newcommand{\R}{\mathbb{R}}

\newcommand{\N}{\mathbb{N}}
\newcommand{\M}{\mathbb{M}}
\newcommand{\cF}{\mathcal{F}}

\newcommand{\bye}{\end{document}}
\newcommand{\by}{\end{proof}\bye}

\newcommand{\hello}{\begin{document}}
\newcommand{\fr}{\frac} 
\newcommand{\disp}{\displaystyle}  
\newcommand{\ga}{\alpha}     
\newcommand{\go}{\omega}
\newcommand{\gep}{\varepsilon}      
\newcommand{\ep}{\gep}    
\renewcommand{\mid}{\,:\,}   
\newcommand{\gb}{\beta} 
\newcommand{\gam}{\gamma}
\newcommand{\gd}{\delta}
\newcommand{\gz}{\zeta} 
\newcommand{\gth}{\theta}   
\newcommand{\gk}{\kappa} 
\newcommand{\gl}{\lambda}
\newcommand{\gL}{\Lambda}
\newcommand{\gs}{\sigma}   
\newcommand{\gf}{\varphi}                  
\newcommand{\tim}{\times}                        
\newcommand{\aln}{&\,}
\newcommand{\ol}{\overline}
\newcommand{\ul}{\underline}           
\newcommand{\pl}{\partial}
\newcommand{\hb}{\text}                
\newcommand{\Int}{\mathop{\text{int}}}

\newcommand{\gG}{\varGamma}
\newcommand{\lan}{\langle}
\newcommand{\ran}{\rangle}
\newcommand{\cD}{\mathcal{D}}
\newcommand{\cB}{\mathcal{B}}
\newcommand{\bcases}{\begin{cases}}
\newcommand{\ecases}{\end{cases}}
\newcommand{\balns}{\begin{align*}}
\newcommand{\ealns}{\end{align*}}
\newcommand{\balnd}{\begin{aligned}}
\newcommand{\ealnd}{\end{aligned}}
\newcommand{\1}{\mathbf{1}}
\newcommand{\bproof}{\begin{proof}}

\newcommand{\eproof}{\end{proof}}

\theoremstyle{definition}
\newtheorem{definition}{Definition}
\theoremstyle{plain}
\newtheorem{theorem}[definition]{Theorem}
\newtheorem{corollary}[definition]{Corollary}
\newtheorem{lemma}[definition]{Lemma}
\newtheorem{proposition}[definition]{Proposition}
\theoremstyle{remark}
\newtheorem{remark}[definition]{Remark}
\newtheorem{notation}[definition]{Notation}

\newcommand{\red}[1]{\textcolor{red}{#1}}
\newcommand{\blu}[1]{\textcolor{blue}{#1}}
\newcommand{\beq}{\begin{equation}}
\newcommand{\eeq}{\end{equation}}
\newcommand{\bthm}{\begin{theorem}}
\newcommand{\ethm}{\end{theorem}}

\usepackage[abbrev]{amsrefs}
\newcommand\coolrightbrace[2]{%
\left.\vphantom{\begin{matrix} #1 \end{matrix}}\right\}#2}
\newcommand{\eqr}[1]{\eqref{#1}}
\newcommand{\bmat}{\begin{pmatrix}}
\newcommand{\emat}{\end{pmatrix}}

\newcommand{\Pmo}{\mathcal{P}^-_{1}}
\newcommand{\Ppo}{\mathcal{P}^+_{1}}
\newcommand{\Pmk}{\mathcal{P}^-_{k}}
\newcommand{\Ppk}{\mathcal{P}^+_{k}}
\newcommand{\diag}{\operatorname{diag}}
\newcommand{\rT}{\mathrm{T}}
\newcommand{\Rn}{{\mathbb R}^N}
\newcommand{\IN}{\text{ in } }
\newcommand{\AND}{\text{ and }}
\newcommand{\FOR}{\text{ for }}
\newcommand{\IF}{\hb{ if }}
\newcommand{\I}{\mathbb{I}}
\newcommand{\du}[1]{\lan#1\ran}
\newcommand{\bald}{\begin{aligned}}
\newcommand{\eald}{\end{aligned}}
\newcommand{\stm}{\setminus}
\newcommand{\lip}{\operatorname{lip}}
\newcommand{\B}{\operatorname{\mathbb{B}}}
\newcommand{\CONV}{\operatorname{co}}

\newcommand{\Sp}{\operatorname{Sp}}
\renewcommand{\P}{\operatorname{\mathbb{P}}}
\newcommand{\erf}{\eqref}
\newcommand{\cG}{\mathcal{G}}
\newcommand{\pr}{\,^\prime}
\newcommand{\gX}{\Xi}
\newcommand{\cV}{\mathcal{V}}
\newcommand{\tH}{\widetilde H}
\renewcommand{\subjclassname}{%
\textup{2010} Mathematics Subject Classification}

\title[Vanishing discount problem]{The vanishing discount problem for \\ 
monotone systems
of Hamilton-Jacobi equations. \\ Part 1: linear coupling}

\author[H. Ishii]{Hitoshi Ishii}
\address[\textsc{Hitoshi Ishii}]{Institute for Mathematics and Computer Science\newline
\indent Tsuda University  \newline
 \indent   2-1-1 Tsuda-machi, Kodaira-shi, Tokyo, 187-8577 Japan.
}
\email{hitoshi.ishii@waseda.jp}

\dedicatory{Dedicated to Italo Capuzzo Dolcetta 
with respect, admiration, and friendship on the occasion of his retirement.}

\begin{document}
\maketitle
\begin{abstract}
We establish a convergence theorem for the vanishing discount problem 
for a weakly coupled system of Hamilton-Jacobi equations. 
The crucial step is the introduction of Mather measures and their relatives for the system, which we call respectively viscosity Mather and Green-Poisson measures. 
This is done by the convex duality and the duality between the space of continuous 
functions on a compact set and the space of Borel measures on it. 
This is part 1 of our study of the vanishing discount problem for systems, which 
focuses on the linear coupling, while part 2 will be concerned with nonlinear coupling.   
\end{abstract}

\section{Introduction}

We consider the weakly coupled $m$-system of Hamilton-Jacobi equations
\[\tag{P$_\gl$}\label{Pl}
\gl v^\gl+Bv^\gl+H[v^\gl]=0 \ \ \IN \T^n, 
\]
where $m\in\N$, $\gl$ is a nonnegative constant, 
called the discount factor in terms of 
optimal control. 
Here $\T^n$ denotes the $n$-dimensional flat torus, $H=(H_i)_{i\in\I}$ is a family of Hamiltonians given by 
\beq\label{H}\tag{H}
H_i(x,p)=\max_{\xi\in\gX}[-g_i(x,\xi)\cdot p -L_i(x,\xi)],
\eeq
where $\I=\{1,\ldots,m\}$, $\gX$ is a given compact metric space, 
$g=(g_i)_{i\in\I}\in C(\T^n\tim\gX,\R^n)^m$ 
and $L=(L_i)_{i\in\I}\in C(\T^n\tim\gX)^m$. 
The unknown in \erf{Pl} is an $\R^m$-valued function $v^\gl=(v^\gl_i)_{i\in\I}$ on 
$\T^n$, $B\mid C(\T^n)^m \to C(\T^n)^m$ 
is a linear  map represented by 
a matrix $B=(b_{ij})_{i,j\in\I}\in C(\T^n)^{m\tim m}$, that is,  
\[
(Bu)_i(x)=(B(x)u(x))_i:=\sum_{j\in\I}b_{ij}(x)u_j(x) \ \ \FOR (x,i)\in \T^n\tim \I.
\]
We use the abbreviated notation $H[v^\gl]$ to denote $(H_i(x,Dv_i^\gl(x))_{i\in\I}$. 
The system is called weakly coupled since the $i$-th equation depends on $Dv^\gl$ only through $Dv_i^\gl$ but not on $Dv_j^\gl$, with $j\not=i$. 
Problem 
\erf{Pl} can be stated in the component-wise style as
\[
\gl v_i^\gl+\sum_{j\in\I}b_{ij}(x)v_j^\gl+H_i(x,Dv_i^\gl)=0 \ \ \IN \T^n,\ i\in\I.
\]

We are mainly concerned with the asymptotic behavior of the solution $v^\gl$ of
\erf{Pl} as $\gl \to 0+$. Asymptotic problems in this class are called the vanishing 
discount problem, in view that the constant $\gl$ in \erf{Pl} appears as a discount factor 
in the dynamic programming PDE in optimal control. 

Recently, there has been a keen interest in the vanishing discount problem 
concerned with Hamilton-Jacobi equations and, furthermore, fully nonlinear degenerate elliptic PDEs. We refer to \cite{DFIZ, AAIY, IJ, MT, IMT1, IMT2, GMT,
CCIZ, DZ2, IS} for relevant work. The asymptotic analysis in these papers relies heavily on Mather measures 
or their generalizations and, thus, it is considered part of the Aubry-Mather and weak KAM theories. For the development of these theories we refer to 
\cite{Fa1, Fa2, Ev1} and the references therein. 

We are here interested in the case of systems of Hamilton-Jacobi equations and, indeed, Davini and Zavidovique in \cite{DZ2} have established 
a general convergence result for the vanishing discount problem for \erf{Pl}. 
We establish a result (Theorem~\ref{thm4.1} below) similar to the main result of \cite{DZ2}.  In establishing our convergence result, we adapt the argument in 
\cite{IMT1} (see also \cite{Go}) to the case of systems, especially, to construct generalized Mather 
measures for \erf{Pl}. Regarding the recent developments of the weak KAM theory and asymptotic analysis 
in its influence for systems, we refer to \cite{CGT,CLLN,MSTY,MT1, MT2, MT3, Ter}. 

The new argument, which is different from that of \cite{DZ2}, makes it fairly easy to build a 
generalized Mather measure for systems in a wide class. 
One advantage of our argument is 
that it allows us to treat the case where the coupling matrix $B$ in \erf{Pl} 
depends on the space variable $x\in\T^n$.  
As in \cites{IJ, IMT1}, our approach is applicable to the system 
with nonlinear coupling of fully nonlinear 
second-order elliptic PDEs, but we restrict ourselves in this paper to the case of the linearly coupled system of first-order 
Hamilton-Jacobi equations.  
Another possible approach for constructing 
generalized Mather measures is the so-called adjoint method (see 
\cite{Ev2, MT, GMT, CGT, Ter}).  

This paper is part 1 of our study of the vanishing discount problem for weakly coupled systems of Hamilton-Jacobi equations and deals only with the linear coupling
and with compact control sets $\Xi$. These restrictions make the presentation 
of our results clear and transparent. In part 2 \cite{IJ}, we remove these restrictions 
and establish a general convergence result extending 
Theorem \ref{thm4.1} below. Sections 5 and 6 are devoted to the study of 
ergodic problems of the form $Bu+H[u]=c$, where $c\in\R^m$ is 
an unknown as well. Also, thanks to the linearity of the coupling,  our results on the ergodic problems are applied to extend the scope of Theorem \ref{thm4.1}.  On the other hand, the role of the ergodic problem, with general right-hand side $c$, is not clear at least for the author 
in the vanishing discount problem for the systems with the nonlinear coupling.

In this paper, we adopt the notion of viscosity solution to \erf{Pl}, for which 
the reader may consult \cite{BaCa,BarB, CIL, PLL}.

To proceed, we give our main assumptions on the system \erf{Pl}. 

We assume that $H$ is coercive, that is, for any $i\in\I$, 
\beq\tag{C}\label{C}
\lim_{|p|\to\infty} \min_{x\in\T^n}H_i(x,p) =\infty. 
\eeq
This is a convenient assumption, under which any upper semicontinuous subsolution of  
\erf{Pl} is Lipschitz continuous on $\T^n$.  

We assume that $B(x)=(b_{ij}(x))$ is a monotone matrix for every $x\in\T^n$, that is, it satisfies
\beq\tag{M}\label{M}
\left\{\begin{minipage}{0.8\textwidth} for any $x\in\T^n$,\  if \
$u=(u_i)_{i\in\I}\in\R^m$ \ and \ $u_k=\max_{i\in\I}u_i\geq 0$, \ then \
$(B(x)u)_k\geq 0$.
\end{minipage} \right.
\eeq
This is a natural assumption that \erf{Pl} 
should possess the comparison principle 
between a subsolution and a supersolution. 
 
In what follows we set, for $\gl\geq 0$,
\[
B^\gl =\gl I+B,
\]
and \erf{Pl} can be written as
\[
B^\gl v^\gl+H[v^\gl]=0 \ \ \IN \T^n.
\]
We use the symbol $u\leq v$ (resp., $u\geq v$) for $m$-vectors $u,v\in\R^n$ to indicate 
$u_i\leq v_i$ (resp., $u_i\geq v_i$) for all $i\in\I$.   

The following theorem is well-known: see \cite{EL, IK} for instance.

\begin{theorem}\label{thm1.1} Assume \eqref{C} and \eqref{M}. Let $\gl>0$. Then the exists a unique 
solution $v^\gl\in\Lip(\T^n)^m$ of \erf{Pl}. Also, if $v=(v_i), w=(w_i)$ are, 
respectively, upper and lower semicontinuous on $\T^n$ and a subsolution and 
a supersolution of \erf{Pl}, then $v\leq w$ on $\T^n$. 
\end{theorem} 

Henceforth, let $\1$ denote the vector $(1,\ldots,1)\in\R^m$. 

\bproof[Outline of proof] We follow the line of the arguments in \cite{IK}. 
Although \cite{IK} is concerned with the case when the domain is an open subset of a Euclidean space,  the results in \cite{IK} is valid in the case when 
the domain is $\T^n$. 

Choose a large constant $C>0$ so that the constant functions $\pm C\1$ 
are a supersolution and a subsolution of \erf{Pl}, respectively. (See also \erf{effect of 1} below.)
According to \cite[Theorems~3.3, Lemma~4.8]{IK}, there is a function $v^\gl=(v_i^\gl)_{i\in\I}
\mid \T^n \to \R^m$ such that the upper and lower semicontinuous envelopes 
$(v^\gl)^*$ and $v^\gl_*$ are a subsolution and a supersolution of \erf{Pl}, 
respectively. By the coercivity assumption \erf{C}, we find (see \cite[Theorem~I.14]{CL}, \cite[Example 1]{IPerron}) 
that the functions $(v_i^\gl)^*$ are Lipschitz continuous on $T^n$.  
Let $R_1>0$ be a Lipschitz bound of the functions $(v_i^\gl)^*$. 
To take into account the Lipschitz property of $(v_i^\gl)^*$, we modify 
the Hamiltonian $H$.  
Fix any $M>0$ so that 
\beq\label{growth}
\max_{(x,\xi,i)\in\T^n\tim\gX\tim\I}|g_i(x,\xi)|<M,
\eeq
and choose constants $N>0$ and $R_2>0$ so that 
\[
H_i(x,p)\geq M|p|-N \ \ \FOR (x,p,i) \in\T^n\tim B_{R_1}\tim\I,
\]
and, in view of \erf{growth},
\[
H_i(x,p)\leq M|p|-N \ \ \FOR (x,p,i) \in\T^n\tim B_{R_2}\tim\I. 
\]
Define $G=(G_i)_{i\in\I}\in C(\T^n\tim\R^n)^m$ by
\[
G_i(x,p)=H_i(x,p)\vee (M|p|-N). 
\]
By the choice of $R_1$, it is easy to see that $(v^\gl)^*$ is a 
subsolution of 
\beq\label{Geq}
\gl u+Bu+G[u]=0 \ \ \T^n.
\eeq
Also, since $G\geq H$, $v_*^\gl$ is a supersolution of \erf{Geq}. 
Observe furthermore that, if $|p|\geq R_2$, then 
\[
G_i(x,p)=M|p|-N \ \ \FOR (x,i)\in\T^n\tim\I,
\] 
the functions $G_i$ are uniformly continuous on $\T^n\tim B_{R_2}$, 
and hence, for some continuous function $\go$ on $[0,\,\infty)$, with $\go(0)=0$, 
\[
|G_i(x,p)-G_i(y,p)|\leq \go(|x-y|) \ \ \FOR (x,y,p)\in(\T^n)^2 \tim\R^n,\,i\in\I.
\] 
The last inequality above shows that $G$ satisfies \cite[(A.2)]{IK}, which allows us to apply \cite[Theorem~4.7]{IK}, to conclude that $(v^\gl)^*\leq v^\gl_*$ on $\T^n$ 
and, moreover, that $v^\gl\in \Lip(\T^n)^*$. Similarly, we deduce that the comparison 
assertion is valid.  Thus, $v^\gl$ is a unique solution of \erf{Pl}. 
\eproof

Regarding the coercivity \erf{C}, the following proposition is well-knwon.

\begin{proposition} The function given by \erf{H} satisfies \erf{C} if and only if 
there exists $\gd>0$ such that 
\beq\label{co}
B_\gd \subset \CONV \{g_i(x,\xi)\mid \xi\in\gX\} \ \ \FOR (x,i)\in\T^n\tim\I,
\eeq
where $\,\CONV$ designates ``convex hull'' and $B_\gd$ denotes the open 
ball with origin at the origin and radius $\gd$.   
\end{proposition}

\bproof[Outline of proof]  Set 
$C(x,i)=\CONV \{g_i(x,\xi)\mid \xi\in\gX\}.$ 
Assume that \erf{co} is valid for some $\gd>0$ and  
observe that
\[\bald
H_i(x,p,u)&\,\geq \max_{\xi\in\gX}-g_i(x,\xi)\cdot p -\max_{(x,i,\xi)\in\T^n\tim\I\tim\gX}L_i(x,\xi)
\\&\,=\max_{q\in C(x,i)} -q\cdot p-\max_{(x,i,\xi)\in\T^n\tim\I\tim\gX}L_i(x,\xi)
\geq \sup_{q\in B_\gd}-q\cdot p- \max_{(x,i,\xi)\in\T^n\tim\I\tim\gX}L_i(x,\xi)
\\&\,=\gd|p|-\max_{(x,i,\xi)\in\T^n\tim\I\tim\gX}L_i(x,\xi),
\eald\]
which shows that $\erf{C}$ holds.  

Next, assume that \erf{co} does not hold for any $\gd>0$. Then there 
exists $(x_k,i_k)\in\T^n\tim\I$ for each $k\in\N$ such that  
\[
B_{1/k}\stm C(x_k,i_k)\not=\emptyset. 
\]
For each $k\in\N$ select $q_k\in B_{1/k}\stm C(x_k,i_k)$ 
and $r_k\in C(x_k,i_k)$ so that $r_k$ is the point of $C(x_k,i_k)$ closest to $q_k$. 
(Notice that $C(x_k,i_k)$ is a compact convex set.) Setting $\nu_k=(q_k-r_k)/|q_k-r_k|$, 
we find that 
\[
\nu_k\cdot (q-r_k)\leq 0 \ \ \FOR q\in C(x_k,i_k).
\]
Sending $k\to\infty$ along an appropriate subsequence, say $(k_j)_{j\in\N}$, 
we find that there are 
a unit vector $\nu=\lim_{j\to\infty} \nu_{k_j}$ of $\R^n$, $r=\lim_{j\to\infty}r_{k_j}\in \R^n$ and $(x,i)\in\T^n\tim\I$ such that
\[
r\in C(x,i) \ \ \AND \ \ 
\nu\cdot (q-r)\leq 0 \ \ \FOR q\in C(x,i). 
\]
If $r\not=0$, then we have $\nu=-r/|r|$, since $\lim_{k\to\infty}q_k=0$, and the inequality above reads
\[
\nu\cdot q\leq -|r|<0 \ \ \FOR q\in C(x,i). 
\]
These observations imply that for $t>0$,
\[
H_{i}(x,-t\nu)=\max_{\xi\in\gX} tg_i(x,\xi)\cdot \nu -\min_{\xi\in\gX} L_i(x,\xi)
\leq -\min_{\xi\in\gX} L_i(x,\xi),
\]
which shows that \erf{C} does not hold. This completes the proof. 
\eproof

The rest of this paper is organized as follows. In Section 2,
we recall some basic facts concerning monotone matrices. In Section 3, 
we study viscosity Green-Poisson measures for our system, which are crucial 
in our asymptotic analysis. We establish the main result for the vanishing discount 
problem in Section 4. We study the ergodic problem (P$_0$) in the cases when $B$ is irreducible, and $B$ is a constant matrix, respectively, in Sections 5 and 6, and combine the results with the analysis on the vanishing discount problem of Section 4.

\section{Monotone matrices} 

Here we are concerned with $m\tim m$ real matrix $B=(b_{ij})_{i,j\in\I}$.

Let $e_i$ denote the vector $(e_{i1},\ldots,e_{im})$, with 
$e_{ii}=1$ and $e_{ij}=0$ if $i\not=j$. 

\begin{lemma} \label{thm2.1} Let $B=(b_{ij})$ be a real $m\tim m$ matrix. 
It is monotone if and only if  
\beq\label{thm2.1.1}
b_{ij}\leq 0 \ \ \IF i\not= j \ \ \AND \ \ \sum_{j\in\I} b_{ij}\geq 0 \ \ \FOR i\in\I. 
\eeq
\end{lemma}

We remark that if $B$ satisfies \erf{thm2.1.1}, then 
\beq \label{thm2.1.2}
b_{ii}=\sum_{j\in\I}b_{ij}-\sum_{j\neq i}b_{ij}\geq 0.
\eeq

\bproof We assume first that $B$ is monotone. 
Since 
\[
\1_{i}=1=\max_{j}\1_{j}>0,
\]
By the monotonicity of $B$, we have
\beq\label{effect of 1}
0\leq (B\1)_{i}=\sum_{j=1}^m b_{ij}\1_j=\sum_{j=1}^m b_{ij} \ \ \FOR i\in\I. 
\eeq
Similarly, if $i\not=j$ and $t\geq 0$, then we have $1=(e_i-te_j)_i=\max_{k\in\I}(e_i-te_j)_k$ and hence,  
\[
0\leq (B(e_i-te_j))_i=b_{ii}-tb_{ij},
\]
from which we find by sending $t\to \infty$ that 
\[
 b_{ij}\leq 0. 
\]
Hence, \erf{thm2.1.1} is satisfied. 

Next, we assume that \erf{thm2.1.1} holds. Let $u\in\R^m$ satisfy
\[
u_k=\max_{i\in\I}u_i\geq 0.
\]
Then we observe that, since $u_k\geq u_j$ for all $j\in\I$,
\[
(Bu)_k=\sum_{j\in\I}b_{kj}u_j=b_{kk}u_k +\sum_{j\neq k}b_{kj}u_j
=b_{kk}u_k+\sum_{j\neq k}b_{kj}u_k=u_k\sum_{j\in\I}b_{kj}\geq 0.
\]
Thus, $B$ is monotone. 
\eproof

\begin{lemma} \label{thm2.2}
Let $u\in\R^m$ and $C\geq 0$ be a constant. Let $B$ be an 
$m\tim m$ 
real monotone matrix. Then we have 
\[
B(u-C\1)\leq Bu\leq B(u+C\1).
\]
\end{lemma}

\bproof  Using Lemma~\ref{thm2.1}, we see that 
\[
(B\1)_i=\sum_{j\in\I}b_{ij}\geq 0 \ \ \FOR i\in\I,
\]
which states that $B\1\geq 0$.  It is then obvious to compute that
\[
B(u+C\1)-Bu=CB\1, \quad Bu-B(u-C\1) =CB\1 \ \ \AND \ \ CB\1\geq 0
\]
and therefore, 
\[
B(u+C\1)\geq Bu\geq B(u-C\1).\qedhere
\]
\eproof

\section{Viscosity Green-Poisson measures}

For $\gl\geq 0$ we write $\cF(\gl)$ for the set of  all $(\phi,u)\in C(\T^n\tim \gX)^m\tim C(\T^n)^m$ such that  $u$ is a subsolution  of  
\[
B^\gl u+H_{\phi}[u]=0 \ \ \IN \T^n,
\]
where $H_{\phi}=(H_{\phi,i})_{i\in\I}$ and 
\[
H_{\phi.i}(x,p)=\max_{\xi\in\gX}(-g_i(x,\xi)\cdot p -\phi_i(x,\xi)).
\]

In the above, since $\phi$ is bounded on $\T^n\tim\gX$, 
if $H$ satisfies \erf{C}, then $H_\phi$ satisfies \erf{C}.

\begin{lemma} \label{cone} The set $\cF(\gl)$ is a convex cone in 
$C(\T^n\tim \gX)^m\tim C(\T^n)^m$ with vertex 
at the origin.
\end{lemma}

\bproof Recall \cite[Remark 2.5]{Bar} that for any $u\in Lip(\T^n)^m$, 
$u$ is a subsolution of 
\[
B^\gl u+H[u]=0 \ \ \IN\T^n
\]
if and only if for any $i\in\I$, 
\[
(B^\gl u)_i(x)+H_i(x,Du_i(x))\leq 0 \ \ \text{ a.e. in }\T^n,
\]
and by the coercivity \erf{C} that for any $(\phi,u)\in\cF(\gl)$, we have
$u\in\Lip(\T^n)^m$.

Fix $(\phi,u), (\psi,v)\in\cF(\gl)$ and $t,s\in[0,\infty)$. Fix $i\in\I$ and 
observe that
 \[\bald
&\,(B^\gl u)_i(x)+H_{\phi,i}(x,Du_i(x))\leq 0 \ \ \text{ a.e. in }\T^n,
\\&\,(B^\gl v)_i(x)+H_{\psi,i}(x,Dv_i(x))\leq 0 \ \ \text{ a.e. in }\T^n,
\eald
\]
which imply that there is a set $N\subset\T^n$ of Lebesgue measure zero such that  
\[\bald
&\,(B^\gl u)_i(x)\leq g(x,\xi)\cdot Du_i(x) +\phi_i(x.\xi)  \ \ \text{ for all } (x,\xi)
\in\T^n\stm N\,\tim\,\gX,
\\&\,(B^\gl v)_i(x)\leq g_i(x,\xi)\cdot Dv_i(x)+\psi_i(x,\xi) \ \ \text{ for all }(x,\xi)\in\T^n\stm N\,\tim\,\gX.
\eald
\]
Multiplying the first and second by $t$ and $s$, respectively, adding 
the resulting inequalities and setting $w=tu+sv$, we obtain
\[
(B^\gl w)_i(x)\leq g(x,\xi)\cdot Dw_i(x) +(t\phi_i+s\psi)(x.\xi)  \ \ \text{ for all } (x,\xi)
\in\T^n\stm N\,\tim\,\gX,
\]
which readily implies that $t(\phi,u)+s(\psi,v)\in\cF(\gl)$. 
\eproof

We refer the reader to \cite[Lemma~2.2]{IMT1} for another proof of the above lemma.

We establish a representation formula for the solution of \erf{Pl}, with $\gl>0$, by 
modifying the argument in \cite{IMT1} (see also \cite{Go}).

For any nonnegative Borel measure $\nu$ on $\T^n\tim\gX$ and 
$\phi\in C(\T^n\tim\gX)$, we write 
\[
\du{\nu,\phi}=\int_{\T^n\tim\gX}\phi(x,\xi)\nu(dx,\,d\xi).
\] 
Similarly, for any collection $\nu=(\nu_i)_{i\in\I}$ of 
nonnegative Borel measures on $\T^n\tim\gX$ and 
$\phi=(\phi_i)\in C(\T^n\tim\gX)^m$, 
we write 
\[
\du{\nu,\phi}=\sum_{i\in\I}\du{\nu_i,\phi_i}\in\R. 
\]
Note that any collection $\nu=(\nu_i)_{i\in\I}$ of 
nonnegative Borel measures on $\T^n\tim\gX$ is regarded as a 
nonnegative Borel measure on $\T^n\tim\gX\tim\I$ and vice versa.

We set 
\[
\rho_i(x):=\sum_{j\in\I}b_{ij}(x) \ \ \FOR i\in\I.
\]
Note that 
\beq\label{thm3.1.0}
B\1=\bmat b_{11}(x)&\cdots& b_{1m}(x)\\ 
\vdots&& \vdots \\
b_{m1}(x)&\cdots&b_{mm}(x)
\emat\bmat 1\\\vdots\\ 1\emat=\bmat \rho_1(x)\\ \vdots\\ \rho_m(x) \emat \ \ \AND \ \ B^\gl\1 =
\bmat \gl+\rho_1(x) \\ \vdots \\ \gl+\rho_m(x)\emat. 
\eeq
By assumption \eqref{M} and Lemma~\ref{thm2.1}, we have 
$\rho_i\geq 0$ on $\T^n$ for all $i\in\I$.

Given a constant $\gl>0$, let $\P_{B^\gl}$ denote the set of  
of nonnegative Borel measures $\nu=(\nu_i)_{i\in\I}$ on 
$\T^n\tim\gX\tim\I$ such that 
\[
\du{\nu,B^\gl\1}=1.
\] 
The last condition reads 
\[
\sum_{i\in\I}(\gl |\nu_i|+\du{\nu_i,\rho_i})=1,
\]
where $|\nu_i|$ denotes the total mass of $\nu_i$ on $\T^n\tim\gX$. 
Note as well that $\P_{B^\gl}$ can be identified with the space of 
Borel probability measures on $\T^n\tim\gX\tim\I$ by the correspondence 
between $\nu=(\nu_i)_{i\in\I}$ and $\sum_{i\in\I}(\gl+\rho_i)\nu_i\otimes \gd_i$, where $\otimes$ indicates the product of two measures and 
$\gd_i$ denotes the Dirac measure at $i$. If we set $\mu:=\sum_{i\in\I}(\gl+\rho_i)\nu_i\otimes \gd_i$ and consider $\mu$ as a collection $(\mu_i)$ of measures 
on $\T^n\tim\gX$, then $\nu_i=(\gl+\rho_i)^{-1}\mu_i$.  
We denote simply by $\P$ the space of Borel probability measures on $\T^n\tim\gX\tim\I$.

For $\gl\geq 0$ and $(z,k)\in\T^n\tim\I$ we set 
\[
\cG(z,k,\gl):=\{\phi-u_k(z) B^\gl \1 \mid (\phi,u)\in \cF(\gl)\}\subset C(\T^n\tim\gX)^m,
\]
and 
\[
\cG\pr(z,k,\gl)=\{\nu=(\nu_i)_{i\in\I}\in\P_{B^\gl} \mid 
\du{\nu,f} \geq 0 \ \ \FOR \ f=(f_i)\in\cG(z,k,\gl)\}. 
\]

\begin{theorem}\label{thm3.1} Assume \erf{H}, \erf{C} and \erf{M}. Let $\gl> 0$ and $(z,k)\in\T^n\tim\I$. 
Let $v^\gl\in C(\T^n\tim\I)$ 
be the unique solution of \erf{Pl}.  
Then 
there exists a $\nu^{z,k,\gl}
=(\nu^{z,k,\gl}_i)_{i\in\I}\in \cG\pr(z,k,\gl)$ such that 
\beq \label{thm3.1.1}
v^\gl_k(z)=\du{\nu^{z,k,\gl},L}.
\eeq
\end{theorem}

We remark that for any $\nu\in\cG\pr(z,k,\gl)$ 
we have $\du{\nu,L}\geq v_k^\gl(z)\du{\nu,B^\gl\1}=v_k^\gl(z)$ and, accordingly, in the theorem above, 
the measures $\nu^{z,k,\gl}$ has the minimizing property:
\beq\label{minGP}
v^\gl_k(z)=\du{\nu^{z,k,\gl},L}=\min_{\nu\in \cG\pr(z,k,\gl)}\du{\nu,L}. 
\eeq
We call any minimizing family $(\nu_i)_{i\in\I}\in\P_{B^\gl}$ of the optimization problem above a viscosity 
Green-Poisson measure for \erf{Pl}.  

\bproof 
Note first that $(L,v^\gl)\in\cF(\gl)$ 
and hence, for any $\nu\in \cG\pr (z,k,\gl)$, 
\beq\label{thm3.1.1+}
 0\leq\du{\nu,L-v_k^\gl(z)B^\gl\1}=\du{\nu,L}-v_k^\gl(z)\du{\nu,B^\gl\1}
=\du{\nu,L}-v_k^\gl(z). 
\eeq

Next,  we show that 
\beq\label{thm3.1.2}
\sup_{(\phi,u)\in\cF(\gl)}\inf_{\nu\in\P_{B^\gl}}
\du{\nu,L-\phi+(u_k(z)-v_k^\gl(z))B^\gl\1}= 0.
\eeq
Note that for $z\in\T^n$,
\[\bald
\sup_{(\phi,u)\in\cF(\gl)}\inf_{\nu\in\P_{B^\gl}}&
\du{\nu,L-\phi+(u_k(z)-v^\gl_k(z))B^\gl\1}
\\&\geq \inf_{\nu\in \P_{B^\gl}}
\du{\nu,L-\phi+(u_k(z)-v_k^\gl(z))B^\gl\1} \Big|_{(\phi,u)=(L,v^\gl)}=0.
\eald\]
Hence, in order to prove \erf{thm3.1.2}, we only need to show that 
\beq\label{thm3.1.3}
\sup_{(\phi,u)\in\cF(\gl)}\inf_{\nu\in \P_{B^\gl}}
\du{\nu,L-\phi+(u_k(z)-v_k^\gl(z))B^\gl\1} \leq 0.
\eeq

We postpone the proof of \erf{thm3.1.3} and, assuming temporarily that 
\erf{thm3.1.2} is valid,  we prove that there exists $\nu\in \cG\pr(z,k,\gl)$ such that
\beq\label{thm3.1.4}
v_k^\gl(z)=\du{\nu,L},
\eeq
which, together with \erf{thm3.1.1+}, completes the proof.

To prove \erf{thm3.1.4}, we observe that $\P_{B^\gl}$ and, by Lemma~\ref{cone}, $\cF(\gl)$ are convex, 
\[
\P_{B^\gl}\ni \nu\mapsto \du{\nu,L-\phi+(u_k(z)-v_k^\gl(z))B^\gl\1}
\]
is convex and continuous, in the topology of weak convergence of measures, 
for any $(\phi,u)\in\cF(\gl)$ and
\[
\cF(\gl)\ni (\phi,u)\mapsto \du{\nu,L-\phi+(u_k(z)-v_k^\gl(z))B^\gl\1}
\]
is concave and continuous for any $\nu\in\P_{B^\gl}$. 
Hence, noting moreover that $\T^n\tim\gX\tim\I$ is a compact set,  we apply the minimax theorem (\cite{Te,Si}), 
to find from \erf{thm3.1.2} that 
\beq \label{thm3.1.5}\bald
0=\sup_{(\phi,u)\in\cF(\gl)}\min_{\nu\in\P_{B^\gl}}&
\du{\nu,L-\phi+(u_k(z)-v_k^\gl(z))B^\gl\1}
\\&=\min_{\nu\in\P_{B^\gl}}\sup_{(\phi,u)\in\cF(\gl)}
\du{\nu,L-\phi+(u_k(z)-v_k^\gl(z))B^\gl\1}.
\eald
\eeq
Observe by using the cone property of $\cF(\gl)$ that 
\[
\sup_{(\phi,u)\in\cF(\gl)}\du{\nu, u_k(z)B^\gl\1 -\phi}
=\bcases
0 & \IF \ \nu\in \cG\pr(z,k,\gl), \\[3pt]
\infty&\IF \ \nu\in \P_{B^\gl} \stm \cG\pr(z,k,\gl). 
\ecases
\]
This and \erf{thm3.1.5} yield
\[\bald
0&\,=\min_{\nu\in\P_{B^\gl}}\sup_{(\phi,u)\in\cF(\gl)}
\du{\nu,L-\phi+(u_k(z)-v_k^\gl(z)B^\gl\1}
%
\\&\,=\min_{\nu\in\cG\pr(z,k,\gl)}
\sup_{(\phi,u)\in\cF(\gl)}\du{\nu,L-v_k^\gl(z)B^\gl\1}
\\&\,=\min_{\nu\in\cG\pr(z,k,\gl)}
\du{\nu,L-v_k^\gl(z)B^\gl\1}=\min_{\nu\in\cG\pr(z,k,\gl)}
(\du{\nu,L}- v_k^\gl(z)\du{\nu,B^\gl\1})
\\&\,=\min_{\nu\in\cG\pr(z,k,\gl)}
\du{\nu,L}-v_k^\gl(z),
\eald
\]
which proves \erf{thm3.1.4}.

It remains to show \erf{thm3.1.3}. For this, we argue by contradiction and thus suppose that \erf{thm3.1.3} 
does not hold. Accordingly, we have  
\[
\sup_{(\phi,u)\in\cF(\gl)}\inf_{\nu\in\P_{B^\gl}}
\du{\nu,L-\phi+(u_k(z)-v_k^\gl(z))B^\gl\1}> \ep
\]
for some $\ep>0$.  We may select $(\phi,u)\in\cF(\gl)$ so that 
\[
\inf_{\nu\in\P_{B^\gl}}\du{\nu,L-\phi+(u_k(z)-v_k^\gl(z))B^\gl\1}
> \ep. 
\]
That is, for any $\nu\in\P_{B^\gl}$, we have
\[
\du{\nu,L-\phi+(u_k(z)-v_k^\gl(z))B^\gl\1}>\ep=\du{\nu,\ep B^\gl\1}. 
\]
Plugging $\nu=(\gl+\rho_i)^{-1}\gd_{(x,\xi,i)}\in \P_{B^\gl}$, with any $(x,\xi,i)\in\T^n\tim\gX\tim \I$, into the above, we find that 
\[
(L_i-\phi_i)(x,\xi)-(v_k^\gl(z)-u_k(z)-\ep)(B^\gl\1)_i>0.
\]
Hence, we have 
\[
\phi(x,\xi)<L(x,\xi)+(u_k(z)-v_k^\gl(z)-\ep)B^\gl\1 
\ \ \FOR (x.\xi)\in\T^n\tim\R^n.
\]
This ensures that $u$ is a subsolution of 
\[
B^\gl u+H[u]=(u_k(z)-v_k^\gl(z)-\ep)B^\gl\1 \ \ \IN \T^n,
\]
which implies that $u-(u_k(z)-v_k^\gl(z)-\ep)\1$ is a subsolution of \erf{Pl}. 
By comparison (Theorem~\ref{thm1.1}), we get 
\[
u(x)-(u_k(z)-v_k^\gl(z)-\ep)\leq v^\gl(x) \ \ \FOR x\in\T^n.
\]
The $k$-th component of the above, evaluated at $x=z$, yields an obvious contradiction.  Thus we conclude that 
\erf{thm3.1.3} holds. 
\eproof

We have the following characterization of $\cG\pr(z,k,\gl)$. 

\begin{proposition} \label{equiv1} Assume 
\erf{H}, \erf{C} and \erf{M} hold. Let $\nu=(\nu_i)_{i\in\I}\in\P_{B^\gl}$ 
and $(z,k,\gl)\in\T^n\tim\I\tim(0,\infty)$. Then 
we have $\nu\in\cG\pr(z,k,\gl)$ if and only if 
\beq\label{equiv1-1}
\sum_{i\in\I}\du{\nu_i, (B^\gl \psi)_i -g_i\cdot D\psi_i}=\psi_k(z) \ \ \FOR \
\psi=(\psi_i)_{i\in\I} \in C^1(\T^n)^m. 
\eeq
\end{proposition}

\bproof Assume first that $\nu\in \cG\pr(z,k,\gl)$. Fix any 
$\psi=(\psi_i)_{i\in\I} \in C^1(\T^n)^m$ and define $\phi=(\phi_i)_{i\in\I}\in C(\T^n\tim\I)^m$ by
\[
\phi_i(x,\xi)=(B^\gl \psi)_i(x) -g_i(x,\xi)\cdot D\psi_i(x). 
\] 
Observe that $u:=\pm\psi$ satisfy, respectively,  
\[
B^\gl u+H_{\pm\phi}[u]=0 \ \ \IN \T^n,
\]
and, hence, 
\[
\pm(\phi,\psi)\in\cF(\gl). 
\] 
Since $\nu\in \cG\pr(z,k,\gl)$, we have 
\[
\pm \psi_k(z)\leq \du{\nu,\pm\phi}=\pm\du{\nu,\phi},
\]
respectively, which shows that \erf{equiv1-1} is valid. 

Now, assume that \erf{equiv1-1} is satisfied. Fix any $(u,\phi)\in \cF(\gl)$.
As noted in the proof of Theorem~\ref{thm1.1}, we have $u\in\Lip(\T^n)$. 
By the standard mollification technique, given a positive constant $\ep>0$, we can approximate $u$ by a smooth 
function $u^\ep$ so that 
\[
\max_{\T^n}|u-u^\ep|<\ep \ \ \AND 
\ \ B^\gl u^\ep+H_\phi[u^\ep]\leq \ep B^\gl\1 \ \ \IN \T^n. 
\]
The last inequality reads 
\[
B^\gl u^\ep_i(x)-g_i(x,\xi)\cdot Du^\ep_i(x)-\phi_i(x,\xi)\leq \ep (B^\gl\1)_i(x)
\ \ \FOR (x,\xi,i)\in\T^n\tim\R^n\tim\I. 
\]
Integrating the above by $\nu_i$, summing up in $i\in\I$ and using \erf{equiv1-1}, we get 
\[
u^\ep_k(z)-\du{\nu,\phi}\leq \ep\du{\nu,B^\gl\1}=\ep.
\] 
Sending $\ep\to 0$ shows that $\nu\in \cG\pr(z,k,\gl)$. 
\eproof

It is convenient to restate the theorem above as follows. 
For $\mu=(\mu_i)_{i\in\I}\in \P$ and $\gl>0$, consider $\nu=(\nu_i)_{i\in\I}\in\P_{B^\gl}$ given by 
\[
\nu_i:=(\gl+\rho_i)^{-1}\mu_i=\fr{1}{(B^\gl\1)_i}\mu_i. 
\]
(Notice by the above definition that $\du{\nu,B^\gl\1}=\du{\mu,\1}=1$.) 
Observe that for $\phi=(\phi_i)_{i\in\I}\in C(\T^n\tim\gX)^m$,
\[
\du{\nu,\phi}=\sum_{i\in\I}\du{\nu_i,\phi_i}=\sum_{i\in\I}\du{\mu_i,
(\gl+\rho_i)^{-1}\phi_i},
\]
and that for any $(z,k)\in\T^n\tim\I$, we have $\nu\in\cG\pr(z,k,\gl)$ if and only if 
\beq\label{thm3.1.6}
\sum_{i\in\I}\du{\mu_i,(\gl+\rho_i)^{-1}\phi_i}\geq u_k(z) \ \ \FOR (\phi,u)\in
\cF(\gl).
\eeq
The condition above is stated in the spirit of Proposition~\ref{equiv1} as 
\[
\sum_{i\in\I}\du{\mu_i,(\gl+\rho_i)^{-1}((B^\gl\psi)_i
-g_i\cdot D\psi_i)}= \psi_k(z) \ \ \FOR \ \psi=(\psi_i)_{i\in\I}\in C^1(\T^n)^m.
\]
We define
\[
\P(z,k,\gl)=\{\mu=(\mu_i)_{i\in\I}\in\P\mid \mu \text{ satisfies \erf{thm3.1.6}}\}.
\]

The following proposition is an immediate consequence of Theorem~\ref{thm3.1}. 

\begin{corollary} \label{thm3.2}Assume \erf{H}, \erf{C} and \erf{M}. Let $\gl> 0$ and $(z,k)\in\T^n\tim\I$. 
Let $v^\gl\in C(\T^n\tim\I)$ 
be the unique solution of \erf{Pl}. 
Then 
there exists a $\mu^{z,k,\gl}
=(\mu^{z,k,\gl}_i)_{i\in\I}\in \P(z,k,\gl)$ such that 
\beq \label{thm3.2.1}
v^\gl_k(z)=\sum_{i\in\I}\du{\mu_i^{z,k,\gl},(\gl+\rho_i)^{-1}L_i}=
\min_{\mu=(\mu_i)_{i\in\I}\in\P(z,k,\gl)}\ \sum_{i\in\I}\du{\mu_i,(\gl+\rho_i)^{-1}L_i}.
\eeq
\end{corollary}

\section{A convergence result for the vanishing discount problem}

We study the asymptotic behavior of the solution $v^\gl$ of \erf{Pl}, with $\gl>0$, 
as $\gl\to 0$.

We make a convenient assumption on the system (P$_0$): 
\[\tag{E} \label{E}\begin{minipage}{0.8\textwidth}
problem (P$_0$)
has a solution $v_0\in \Lip(\T^n)$. 
\end{minipage}
\]

If $\rho_i>0$ for all $i\in\I$, then Theorem~\ref{thm1.1} assures that there exists 
a \emph{unique} solution $v_0$ of \erf{E}. In this situation, it is not difficult to show that 
the uniform convergence, as $\gl \to 0+$, of $v^\gl$ to the unique solution $v_0$ on $\T^n$. In general, existence and uniqueness of a solution of (P$_0$) may fail. 
In fact, one can prove at least  in the case when the $b_{ij}$ are constants (see Theorem~\ref{thm5.1}) that there exists $c\in\R^m$ such that 
\beq
\label{thm4.0.1} 
Bu+H[u]=c \ \ \IN\T^n
\eeq
has a solution $v_0\in \Lip(\T^n)$ and possibly multiple solutions. 
If such a $c=(c_i)$ exists,  then the introduction of a new family of 
Hamiltonians, 
\[ 
\widetilde H=(\widetilde H_i)_{i\in\I}, \quad\hb{ with }\widetilde H_i(x,p)=H_i(x,p)-c_i,
\]
allows us to view \erf{thm4.0.1} as in the form of (P$_0$). The link between two
vanishing discount problems for the original \erf{Pl} and for \erf{Pl}, with $\tH$ 
in place of $H$, is discussed in Sections 5 and 6.

\begin{theorem} \label{thm4.1}
Assume \erf{H}, \erf{C}, \erf{M} and \erf{E}. Let $v^\gl$ be the unique solution of \erf{Pl} for $\gl>0$. Then there exists a solution $v^0\in\Lip(\T^n)^m$ of \emph{(P$_0$)}
such that  the functions $v_i^\gl$ converge to $v_i^0$ uniformly on $\T^n$ as 
$\gl\to 0$ for all $i\in\I$.
\end{theorem}


\begin{lemma}\label{thm4.2} Under the hypotheses of Theorem~\ref{thm4.1}, 
there exists a constant $C_0>0$ such that for any $\gl>0$,
\beq\label{thm4.2.1}
|v_i^\gl(x)|\leq C_0 \ \ \FOR (x,i)\in\T^n\tim\I.
\eeq
\end{lemma}

\bproof Let $v_0=(v_{0,i})_{i\in\I}\in \Lip(\T^n)^m$ be the solution of (P$_0$). Choose 
a constant $C_0>0$ so that 
\[
|v_{0,i}(x)|\leq C_1 \ \ \FOR (x,i)\in\T^n\tim\I,
\]
and observe by the monotonicity of $B$ (Lemma~\ref{thm2.2}) that $v_0+C_1\1$ 
and $v_0-C_1\1$ are a supersolution and a subsolution of (P$_0$), respectively. 
Noting that $v_0+C_1\1\geq 0$ and $v_0-C_1\1\leq 0$, we deduce that 
$v_0+C_1\1\geq 0$ and $v_0-C_1\1\leq 0$ are a supersolution and a subsolution of 
\erf{Pl} for any $\gl>0$, respectively. By comaprison (Theorem \ref{thm1.1}), 
we see that, for any $\gl>0$, 
$v_0-C_1\1\leq v^\gl\leq v_0+C_1\1$ on $\T^n$ and, moreover, 
$-2C_1\1\leq v^\gl\leq 2C_1\1$ on $\T^n$. Thus, \erf{thm4.2.1} holds with 
$C_0=2C_1$.
\eproof 

\begin{lemma}\label{thm4.3} Under the hypotheses of Theorem~\ref{thm4.1}, 
the family $(v^\gl)_{\gl\in(0,\,1)}$ is equi-Lipschitz continuous on $\T^n$. 
\end{lemma}

Indeed, the family $(v^\gl)_{\gl>0}$ is equi-Lipschitz continuous on $\T^n$, which  
we do not need here. 

\bproof By Lemma~\ref{thm4.2}, there is a constant $C_0>0$ such that 
\[
|(B^\gl v^\gl(x))_i|\leq C_0 \ \ \FOR (x,i,\gl)\in\T^n\tim\I\tim(0,\,1). 
\]
Hence, as $v^\gl$ is a solution of \erf{Pl},   we deduce by \erf{C} that there 
exists a constant $C_1>0$ such that the $v_i^\gl$ are subsolutions of 
$|Du|\leq C_1$ in $\T^n$. It is a standard  fact that the $v_i^\gl$ are 
Lipschitz continuous on $\T^n$ with $C_1$ as their Lipschitz bound.  
\eproof

In the proof of Theorem~\ref{thm4.1}, Corollary~\ref{thm3.2} has a crucial role. 
We need also results for $\gl=0$ similar to the corollary.

We consider the condition for $\mu\in\P$, 
\beq\label{thm3.2.2}
\du{\mu,\phi}\geq 0 \ \ \FOR (\phi,u)\in\cF(0).
\eeq
We denote by $\P(0)$ the subset of $\P$ consisting of those $\mu$ which satisfy
\erf{thm3.2.2}.

\begin{theorem}\label{thm3.3}
Assume \erf{H}, \erf{C}, \erf{M} and \erf{E}. Assume that $\rho_i=0$ on $\T^n$ for every $i\in\I$. Then there 
exists a $\mu^0=(\mu_i^0)_{i\in\I}\in\P(0)$ such that
\beq\label{thm3.3.1}
0=\du{\mu^0,L}=\min_{\mu\in\P(0)}\du{\mu,L}.
\eeq
\end{theorem}

\bproof We fix a $(z,k)\in \T^n\tim\I$. By Corollary \ref{thm3.2}, for each $\gl>0$
there exists $\mu^\gl=(\mu_i^\gl)_{i\in\I}\in \P(z,k,\gl)$ such that 
\beq\label{thm3.3.2}
\gl v_k^\gl(z)=\sum_{i\in\I} \gl \du{\mu^\gl_i,\gl^{-1}L_i}=\du{\mu^\gl,L}. 
\eeq
Since $(\mu^\gl)_{\gl>0}$ is a family of Borel probability measures on a compact space $\T^n\tim\gX\tim\I$, there exists a sequence $(\gl_j)_{j\in\N}\subset (0,\,1)$ converging to zero such that the sequence $(\mu^{\gl_j})_{j\in\N}$ converges 
weakly in the sense of measures 
to a Borel probability measure $\mu^0$ on $\T^n\tim\gX\tim\I$. 
It follows from \erf{thm3.3.2} and Lemma~\ref{thm4.2} that 
\[
0=\du{\mu^0,L}. 
\]
Observe that if $(\phi,u)\in\cF(0)$, then, for any $\gl>0$, $u$ is a subsolution of 
\[
B^\gl u+H_\phi[u]=\gl u \ \ \IN \T^n,
\] 
and hence, $(\psi,u) \in \cF(\gl)$, with $\psi(x,\xi)=\phi(x,\xi)+\gl u(x)$.  
Hence, the inclusion $\mu^\gl\in\cG\pr(z,k,\gl)$ yields
\[
u_k(z)\leq \sum_{i\in\I}\du{\mu_i^\gl,\gl^{-1}(\phi_i+\gl u_i)}
=\gl\du{\mu^\gl,\phi}+\du{\mu^\gl,u}.
\]
Multiplying the above by $\gl$ and sending $\gl=\gl_j \to 0$, 
in view of Lemma~\ref{thm4.2}, we get 
\[
0\leq \du{\mu^0,\phi}.
\]
This shows that $\mu^0\in\P(0)$. 
These observations together with \erf{thm3.2.2} for $\mu\in\P(0)$ 
guarantee that 
\[
0=\du{\mu^0,L}=\min_{\mu\in\P(0)}\du{\mu,L}. \qedhere 
\] 
\eproof 

We state a characterization of $\P(0)$ in the next, similar to Proposition~\ref{equiv1},
which we leave to the reader to verify.   

\begin{proposition}\label{equiv2} Assume \erf{H}, \erf{C} and \erf{M}. Let $\mu=(\mu_i)_{i\in\I}\in\P$.  We have $\mu\in\P(0)$ if and only if
\[
\sum_{i\in\I}\du{\mu_i,(B\psi)_i-g_i\cdot D\psi_i}=0 \ \ \FOR \
 \psi=(\psi_i)_{i\in\I}\in C^1(\T^n)^m. 
\]
\end{proposition}

We call any minimizer $\mu\in\P(0)$ of the optimization problem \erf{thm3.3.1} a viscosity Mather measure.

We denote by $\M_+$ the set of all Borel nonnegative measures $\mu=(\mu_i)_{i\in\I}$ on $\T^n\tim\gX\tim\I$. We set
\[
\M_+(0)=\{\mu\in\M_+\mid \mu \text{ satisfies \erf{thm3.2.2}}\}.
\] 

\begin{theorem}\label{thm3.4} Let $(z,k)\in\T^n\tim\I$. Assume \erf{H}, 
\erf{C}, \erf{M} and \erf{E}. 
For any $\gl>0$, let $v^{\gl}$ be the unique solution of \erf{Pl} and $\mu^\gl \in\P(z,k,\gl)$ be a minimizer of \erf{thm3.2.1}. Then there exists a subsequence of $(\gl_j)$, which is denoted again by the same symbol, such that, as $j\to\infty$, 
\[
\fr{\gl_j}{\gl_j+\rho_i}\mu_i^{\gl_j} \to \mu_i^0 
\]
weakly in the sense of measures for some $\mu^0=(\mu^0_i)_{i\in\I}\in\M_+(0)$,  
and $\mu^0$ satisfies  
\beq\label{mini0}
\du{\mu^0,L}=0. 
\eeq
In particular, 
\beq\label{thm3.4.1}
0=\du{\mu^0,L}=\min_{\mu\in\M_+(0)}\du{\mu,L}.
\eeq
\end{theorem}

Notice that the minimization problem \erf{thm3.4.1} is trivial since $\mu^0=0$ 
is a minimizer.

\bproof The proof is similar to that of Theorem~\ref{thm3.3}. 

We fix a $(z,k)\in \T^n\tim\I$. For each $\gl>0$, we have 
\beq\label{thm3.4.2}
\gl v_k^\gl(z)=\sum_{i\in\I} \gl \du{\mu^\gl_i,(\gl+\rho_i)^{-1}L_i}. 
\eeq
Observe that 
\[
\du{\gl(\gl+\rho_i)^{-1}\mu^\gl_i,\1}\leq \du{\mu^\gl_i,\1}=\sum_{i\in\I}|\mu_i^\gl|=1.
\]
Accordingly,  since $\T^n\tim\gX\tim\I$ is a compact metric space,
the families $(\gl(\gl+\rho_i)^{-1}\mu_i^\gl)_{\gl=\gl_j, j\in\N}$ 
have a common subsequence, along which all the families converge 
to some Borel nonnegative measures $\mu^0_i$ weakly in the sense of measures. We may assume by 
replacing the original sequence $(\gl_j)$ by its subsequence that 
\[
\fr{\gl_j}{\gl_j+\rho_i}\mu_i^{\gl_j} \to \mu_i^0 
\]
weakly in the sense of measures. Combine this with \erf{thm3.4.2} yields
\[
0=\sum_{i\in\I}\du{\mu^0_i,L_i}=\du{\mu^0,L}. 
\]
It is obvious to see that $\mu^0\in\M_+$. 

Let $(\phi,u)\in\cF(0)$. As before, we have $(\psi,u) \in \cF(\gl)$, with $\psi(x,\xi)=\phi(x,\xi)+\gl u(x)$ and moreover   
\[
u_k(z)\leq \sum_{i\in\I}\du{\mu_i^\gl,(\gl+\rho_i)^{-1}(\phi_i+\gl u_i)}
=\du{\mu^\gl,(\gl+\rho_i)^{-1}\phi}+\gl\du{\mu^\gl,(\gl+\rho_i)^{-1}u}.
\]
Multiplying the above by $\gl$ and sending $\gl=\gl_j \to 0$, we get 
\[
0\leq \du{\mu^0,\phi}.
\]
This shows that $\mu^0\in\M_+(0)$. 
\eproof

\bproof[Proof of Theorem~\ref{thm4.1}] 
Let $\cV$ denote the set of accumulation points $v
=(v_i)\in C(\T^n)^m$ in the space $C(\T^n)^m$ of $v^\gl$ as $\gl\to 0$. 
In view of the Ascoli-Arzela theorem, Lemmas~\ref{thm4.2} and \ref{thm4.3} guarantee that the family $(v^\gl)_{\gl\in(0,\,1)}$ is relatively compact in $C(\T^n)^m$.  In particular, the set $\cV$ is nonempty.  Note by the stability of the viscosity property under uniform 
convergence that any $v\in\cV$ is a solution of (P$_0$). 

If $\cV$ is a singleton, then it is obvious that the whole family 
$(v^\gl)_{\gl>0}$ converges to the unique element of $\cV$ in $C(\T^n)^m$
as $\gl\to 0$. 

We need only to show that $\cV$ is a singleton.  It is enough to show that 
for any $v,w\in\cV$ and $(z,k)\in\T^n\tim\I$, the inequality $w_k(z)\leq v_k(z)$ 
holds. 

Fix any $v,w\in\cV$ and $(z,k)\in\T^n\tim\I$. 
Select sequences $(\gl_j)$ and $(\gd_j)$ converging to zero 
so that
\[
v^{\gl_j} \to v,\ \ v^{\gd_j} \to w \ \ \IN C(T^n)^m \ \ \hb{ as } j\to\infty. 
\] 
By Corollary \ref{thm3.2}, there exists a sequence $(\mu^j)_{j\in\N}$  such that   
\beq\label{thm4.1.1}
\mu^j\in \cG\pr(z,k,\gl_j) \ \ \AND \ \ 
v_k^{\gl_j}(z)=\sum_{i\in\I}\du{\mu^j_i,(\gl_j+\rho_i)^{-1}L_i} \ \ \FOR j\in\N. 
\eeq
In view of Theorem~\ref{thm3.4}, we may assume by passing to a subsequence if necessary that, as $j\to\infty$, 
\[
\fr{\gl_j}{\gl_j+\rho_i}\mu_i^j \to \mu_i^0 \ \ \hb{ weakly in the sense of measures}
\]
for all $i\in\I$ and for some $\mu^0=(\mu^0_i)_{i\in\I}\in\M_+(0)$ and, moreover, 
\beq\label{thm4.1.2}
0=\du{\mu^0,L}. 
\eeq

Since $(L-\gl v^\gl,v^\gl)\in\cF(0)$ and $\mu^0\in\M_+(0)$, in view of \erf{thm4.1.2}, we have 
\[
0\leq \du{\mu^0,L-\gl v^\gl}=
\du{\mu^0,L}-\du{\mu^0, \gl v^\gl}
=-\gl\du{\mu^0,v^\gl},
\]
which yields after dividing  by $\gl>0$ and then sending $\gl \to 0$ 
along $\gl=\gd_j$
\beq\label{thm4.1.3}
\du{\mu^0,w}\leq 0.
\eeq

Now, note that $w$ is a solution of 
\[
B^\gl w+H[w]=\gl w \ \ \IN\T^n,
\]  
and thus, $(L+\gl w, w)\in\cF(\gl)$ and infer by \erf{thm4.1.1} that 
\[
w_k(z)\leq \sum_{i\in\I}\du{\mu^j_i,(\gl_j+\rho_i)^{-1}(L_i+\gl_j w_i)}
=v_k^{\gl_j}(z)+\gl_j\sum_{i\in\I}\du{\mu^j_i,(\gl_j+\rho_i)^{-1}w_i}. 
\]  
Sending $j\to\infty$ now yields
\[
w_k(z)\leq v_k(z)+\du{\mu^0,w}.
\]
This together with \erf{thm4.1.3} shows that $w_k(z)\leq v_k(z)$, which completes 
the proof. 
\eproof

\section{The ergodic problem for irreducible matrix $B$} 

We consider the problem of finding $c=(c_i)_{i\in\I}\in\R^m$ and  
$v=(v_i)_{i\in\I}\in C(\T^n)^m$ such that $v$ is a solution of 
\beq\label{ergodic}
Bv+H[v]=c \ \ \IN \T^n.
\eeq
The pair of such $c$ and $v$ is also called a solution of \erf{ergodic}. 
This problem is called the \emph{ergodic problem} in this paper although 
the term, ergodic problem, should be used only when 
the condition that $\sum_{j\in\I}b_{ij}(x)=0$ holds for some $(i,x)\in\I\tim\T^n$.  

Henceforth, $D(x)$ denotes the diagonal matrix 
\[
D(x)= \diag(\rho_1(x),\ldots,\rho_m(x)) \ \ \ \FOR x\in\T^n,
\]
where, as before, $\rho_i(x)=\sum_{j\in\I}b_{ij}(x)$. 

Throughout this section, we treat the case when
\beq \label{irreducible}
B(x)\ \text{ is irreducible}.
\eeq
The irreducibility of $B(x)$ is stated as follows: for any nonempty subset $I$ of $\I$, which is not identical to $\I$, there exists a pair of $i\in I$ and $j\in\I\stm I$
such that $b_{ij}(x)\not=0$.  

The following result has been established in Davini-Zavidovique
\cite[Theorem 2.10]{DZ1} (see also \cites{CLLN, MT3}).

\begin{proposition}\label{existence-ir} Assume \erf{H}, \erf{C}, \erf{M},  \erf{irreducible}, and that 
\beq \label{degenerate}
\sum_{j\in\I} b_{ij}(x)=0 \ \ \text{ for all }(i,x)\in\I\tim\T^n. 
\eeq
Then there exist $c_0\in\R$ and 
$v_0\in \Lip(\T^n)^m$ such that the pair $(c_0\1,v_0)$ is a solution of \erf{ergodic}. 
\end{proposition}

We remark that \erf{degenerate} is satisfied if and only if $B(x)\1=0$ for all $x\in\T^n$, 
which holds if and only if $\rho_i(x)=0$ for all $(i,x)\in\I\tim\T^n$.

The next theorem states the central result of this section.

\begin{theorem} \label{convergence-ir} Assume \erf{H}, \erf{C}, \erf{M},  
\erf{irreducible}, and \erf{degenerate}. 
Let $v^\gl$ be the unique solution of \erf{Pl} for $\gl>0$. Then there exists a constant $c^0\in\R$ and a function $v^0\in\Lip(\T^n)^m$ such that 
the functions $v^\gl+\gl^{-1}c^0\1$ converge to $v^0$ uniformly on $\T^n$ as 
$\gl\to 0$. Moreover, the pair $(c^0\1,v^0)$ is a solution of \erf{ergodic}. 
\end{theorem}

\bproof Thanks to Proposition \ref{existence-ir}, there exists 
a solution $(c_0,v_0)\in\R^m\tim C(\T^n)^m$ of \erf{ergodic}.  
We set $
\tH=H-c_0\1$, 
and note that, 
since $B(x)\1=0$ for all $x\in\T^n$,
the function $w^\gl:=v^\gl+\gl^{-1}c_0\1$ satisfies, in the viscosity sense, 
\[
\gl w^\gl+Bw^\gl+\tH[w^\gl] 
=\gl v^\gl +c_0\1+Bv^\gl+H[v^\gl]-c_0\1=0. 
\]
By Theorem \ref{thm4.1},  there exists a solution $v^0\in\Lip(\T^n)^m$ of 
$Bv^0+\tH[v^0]=0$ in $\T^n$ such that, as $\gl\to 0+$, 
$w^\gl\to v^0$ in $C(\T^n)^m$.  
Noting that $(c_0\1,v^0)$ is a solution of \erf{ergodic}, we finish
the proof.  
\eproof

The condition \erf{degenerate} in Proposition \ref{existence-ir} can be removed 
and the following  theorem is valid. 

\begin{theorem}\label{existence-ir-2} Assume \erf{H}, \erf{C}, \erf{M}, and \erf{irreducible}. 
Then there exist $c^0\in\R$ and 
$v^0=(v^0_i)_{i\in\I}\in \Lip(\T^n)^m$ such that the pair $(c^0\1,v^0)$ is a solution of \erf{ergodic}. 
\end{theorem}

\bproof For $x\in\I\tim\T^n$, we set 
\[
B^0(x)=(b^0_{ij}(x)):=B(x)-D(x). 
\]
and note that $B^0(x)$ is irreducible and \erf{degenerate} holds with 
$b_{ij}(x)$ replaced by $b^0_{ij}(x)$. Note also that $\rho_i(x)\geq 0$ for all 
$(i,x)\in\I\tim\T^n$. 

Thanks to Proposition \ref{existence-ir}, there exist $c^0\in\R$ and $v=(v_i)\in\Lip(\T^n)^m$ which solve 
\[
B^0v+H[v]=c^0\1 \ \ \IN \T^n. 
\]
We choose a constant $C>0$ so that $\max_{(i,x)\in\I\tim\T^n}|v_i(x)|\leq C$ 
and set $v^\pm(x)=v(x)\pm C\1$, respectively. Observe that, since 
$v^+_i(x)\geq 0$ and $v^-_i(x)\leq 0$ for all $(i,x)\in\I\tim\T^n$, the functions 
$u=v^+$ and $u=v^-$ are a supersolution and subsolution of 
\[
B^0u+Pu+H[u]=c^0\1 \ \ \IN\T^n,
\] 
that is, $Bu+H[u]=c^0\1 \ \ \IN \T^n$, respectively. In view of the Perron method, 
the function $v^0=(v^0_i)_{i\in\I}\in \Lip(\T^n)$ given by 
\[\bald
v^0_i(x)=\sup\{u_i(x)\mid u=(u_i)\in C(\T^n)^m & \text{ is a subsolution of }
Bu+H[u]=c^0\1 \ \IN \T^n, 
\\ & v^-\leq u\leq v^+ \ \IN \T^n\},
\eald
\]
is a solution of \erf{ergodic}, with $c=c^0\1$.
\eproof 

Even without the assumption \erf{degenerate}, it is immediate from Theorem \ref{thm4.1} that, under the hypotheses of Theorem \ref{existence-ir-2}, if $c^0=0$, then the convergence holds  for the whole family of the solutions $v^\gl$ of \erf{Pl}, with $\gl>0$.  A typical case when $c^0=0$ is realized is given by \cite[Theorem 4.2]{CLLN}(see also \cites{DZ1, MT1}). 

\section{The ergodic problem for constant matrix $B$}

Throughout this section we assume that $B$ is a \emph{constant} matrix, that is,  independent of $x\in\T^n$. 

The main results in this section are as follows. 

\begin{theorem}\label{thm5.1}Assume \erf{H}, \erf{C}, \erf{M}, and that $B$ is a 
constant matrix. Then \erf{ergodic} has a solution $(c,v)\in \R^m\tim C(\T^n)^m$. 
\end{theorem}

\begin{theorem} \label{thm5.1.0} Under the same hypotheses of Theorem \ref{thm5.1}, let 
$(c,v_0)\in \R^m\tim C(\T^n)^m$ be a solution of \erf{ergodic} and let $v^\gl$ be the unique solution of \erf{Pl} for $\gl>0$. Then there exists a function $v^0\in C(\T^n)^m$ 
such that  the functions $v^\gl+(\gl I+B)^{-1}c$ converge to $v^0$ uniformly on $\T^n$ as $\gl\to 0$. Moreover, the pair $(c,v^0)$ is a solution of \erf{ergodic}.  
\end{theorem}

\bproof It is well-known (and easily checked) that due to the monotonicity of $B$, 
$(\gl I+B)$ is invertible for any $\gl>0$.  We set $\tH(x,p)=
H(x,p)-c$ \ for $(x,p)\in\T^n\tim\R^n$ and  also
$w^\gl(x)=v^\gl(x)+(\gl I+B)^{-1}c$ \ for $x\in\T^n$.
Observe that, in the viscosity sense,
\[\bald
\gl w^\gl(x)&+Bw^\gl(x)+\tH[w^\gl] 
\\&=\gl v^\gl +Bv^\gl +H[v^\gl]-c+\gl (\gl I+B)^{-1}c+B(\gl I+B)^{-1}c
=0 \ \ \IN \T^n.
\eald\]
It is clear that $\tH$ satisfies \erf{H} and \erf{C} and that $v_0$ is  a solution of 
$Bu+\tH[u]=0$ in $\T^n$.  By Theorem \ref{thm4.1}, we conclude that 
there exists a solution $v^0\in C(\T^n)^m$ of $Bu+\tH[u]=0$ in $\T^n$ 
such that $w^\gl \to v^0$ in $C(\T^n)^m$ as $\gl\to 0+$.  Noting that 
$(c,v^0)$ is a solution of \erf{ergodic}, we finish the proof. 
\eproof

For the proof of Theorem \ref{thm5.1}, we begin with a preliminary remark on 
the permutations. 

For a given permutation $\pi\mid \I \to \I$, 
we define the $m\tim m$ matrix $P$ by
\beq\label{permP}
P=(\gd_{\pi(i),j})_{i,j\in\I},
\eeq
where $\gd_{ij}=\gd_{i,j}:=1$ if $i=j$ and $=0$ otherwise. Note that 
$P^{-1}=(\gd_{i,\pi(j)})_{i,j\in\I}=P^\rT$
and that
for any $u=(u_i)_{i\in\I}$,
\[
Pu=P\bmat u_1\\ \vdots \\ u_m\emat
=\bmat u_{\pi(1)}\\ \vdots\\ u_{\pi(m)}\emat.
\]
The system of Hamilton-Jacobi equations 
\beq\label{sys}
\gl u+Bu+H[u]=0 
\eeq
can be written component-wise as 
\[
\gl u_{\pi(i)}+(Bu)_{\pi(i)}+H_{\pi(i)}[u_{\pi(i)}]
=0 \ \ \FOR i\in\I.
\]
By the use of $P$, the system above is expressed as
\[
\gl (Pu)_i+(PBu)_i+(PH)_{i}[(Pu)_i]=0,
\]
and furthermore, if $v=Pu$,
\beq\label{sysP}
\gl (v)_i+(PBP^\rT v)_i+(PH)_{i}[v_i]=0. 
\eeq
Set $A=(a_{ij})_{i,j\in \I}=PBP^\rT$ and observe that if $B$ is monotone, then
\[
a_{ij}=\sum_{k,l\in\I}\gd_{i,\pi(k)}b_{kl}\gd_{\pi(l),j}=b_{\pi^{-1}(i),\pi^{-1}(j)}
\bcases
\geq 0&\IF i=j,\\
\leq 0&\IF i\neq j,
\ecases 
\]
and 
\[
\sum_{j\in\I}a_{ij}=\sum_{j\in\I}b_{\pi^{-1}(i),\pi^{-1}(j)}=\sum_{j\in\I}b_{\pi^{-1}(i),j}\geq 0.
\]
Consequently, if $B$ is monotone, then $PBP^\rT$ is monotone as well,
and the system \erf{sys}, by using the permutation matrix $P$, is converted 
to \erf{sysP}. 

\bproof[Proof of Theorem \ref{thm5.1}]
It is well-known (see for instance \cite[Section 2.3]{Var}) that, given a monotone matrix $B$,  one can find a permutation $\pi\mid\I\to\I$ such that 
\beq\label{normal}
PBP^\rT=\bmat 
B^{(1)} &0 & \cdots &0 \\
*& B^{(2)}& \ddots& \vdots\\
\vdots &\ddots&\ddots&0\\
*&\cdots &*&B^{(r_p)} \\ 
\emat,
\eeq
where, $P$ is given by \erf{permP}, $B^{(1)}$ is a diagonal matrix of order $r_1$ and, for $1<i\leq p$, 
$B^{(i)}$ are irreducible matrices of order $r_i$.  In view of the 
preliminary remark before this proof, 
to seek for a solution  of \erf{ergodic}, we may 
and do assume henceforth  
$B$ has the normal form of the right hand side of \erf{normal}. 

Set 
\[
s_k=\sum_{1\leq i<k} r_i \ \ \AND \ \ 
\I_k=\{s_k+1,\ldots,s_k+r_k\} \ \ \FOR k\in\{1,\ldots,p\}. 
\]
Notice that $s_1=0$.  
If $r_1\geq 1$, then we first show that there exist an $r_1$-vector 
$c^{(1)}=(c^{(1)}_i)_{i\in\I_1}\in\R^{r_1}$
and a function $v^{(1)}=(v^{(1)}_i)_{i\in\I_1}\in C(\T^n)^{r_1}$ such that 
$v^{(1)}$ is a solution of 
\beq\label{ergo1}
B^{(1)}v^{(1)}+H^{(1)}[v^{(1)}]=c^{(1)} \ \ \IN\T^n,
\eeq
where $H^{(1)}=(H_i)_{i\in \I_1}$. The system is, in fact, a collection of 
single equations 
\beq\label{single}
b_{ii}v_i^{(1)}+H^{(1)}_i[v_i^{(1)}]=c_i^{(1)} \ \ \IN \T^n, \text{ with }i\in\I_1,
\eeq
and thus the existence of a solution  $(c^{(1)},v^{(1)})$  
of \erf{ergo1} is a classical result. Indeed, for each $i\in\I_1$, if $b_{ii}^{(1)}>0$, then  
\erf{single} has a (unique) 
solution $v_i^{(1)}\in \Lip(\T^n)$ for any choice of $c_i^{(1)}$. 
If $b_{ii}^{(1)}=0$, then \erf{single} has a solution $(c_i^{(1)},v_i^{(1)})\in\R\tim \Lip(\T^n)$ (see \cite{LPV}). 
If $r_1=m$, then we are done. 

Next, assume that $r_1<m$ (and equivalently, $1<p$) 
and we show that there exist  a vector 
$c^{(2)}=(c^{(2)}_i)_{i\in\I_2}\in\R^{r_2}$  
and a function $v^{(2)}=(v^{(2)}_i)_{i\in\I_2}\in C(\T^n)^{r_2}$ such that 
$v^{(2)}$ is a solution of 
the system
\beq\label{ergo2}
B^{(2)}v^{(2)}+H^{(2)}[v^{(2)}]=c^{(2)} \ \ \IN \T^n,
\eeq
where 
\beq\label{erg2+1}
H^{(2)}_i(x,p)=H_{i}(x,p)-\sum_{j\in\I_1}b_{i,j}v^{(1)}_j(x) \ \ \FOR i\in\I_2. 
\eeq
According to Proposition~\ref{existence-ir}, there exist 
$c^{(2)}=(c^{(2)}_i)_{i\in\I_2}\in\R^{r_2}$ 
and $v^{(2)}=(v^{(2)}_i)_{i\in\I_2}\in C(\T^n)^{r_2}$ which satisfy \erf{ergo2}.
This way (by induction), we find $c^{(1)},\ldots,c^{(p)}$ and 
$v^{(1)},\ldots,v^{(p)}$ such that 
\[
c^{(k)}\in \R^{r_k} \ \ \AND \ \ v^{(k)}\in C(\T^n)^{r_k} \ \ \FOR k\in\{1,\ldots,p\},
\]
and $v^{(k)}$ satisfies 
\beq\label{ergok}
B^{(k)}v^{(k)}+H^{(k)}[v^{(p)}]=c^{(k)} \ \ \IN \T^n,\ \ \FOR k\in\{1,\ldots,p\}. 
\eeq
where
\beq\label{ergk+1}
H_i^{(k)}(x,p)=H_{i}(x,p)-\sum_{1\leq j<k}\sum_{q\in\I_j} b_{i,q}v^{(j)}_{q}(x) 
\ \ \FOR i\in\I_k. 
\eeq
We define $c=(c_i)_{i\in\I}\in\R^m$ and $v=(v_i)_{i\in\I}\in C(\T^n)^m$ by setting 
\[
c_i=c_i^{(k)} \ \ \AND \ \ v_i=v_i^{(k)} \ \ \FOR i\in\I_k,\, k\in\{1,\ldots,p\}, 
\]
and observe that 
\[
Bv+H[v]=c \ \ \IN \T^n. 
\]
This completes the proof. 
\eproof

\section*{Acknowledgements}
The author would like to thank the anonymous referee for useful and critical comments on the original version of this paper, which have helped significantly to improve the presentation.  This work was partially supported by the KAKENHI \#16H03948, \#18H00833, \#20H03688, JSPS.


\begin{bibdiv}
\begin{biblist}
\bib{AAIY}{article}{
   author={Al-Aidarous, Eman S.},
   author={Alzahrani, Ebraheem O.},
   author={Ishii, Hitoshi},
   author={Younas, Arshad M. M.},
   title={A convergence result for the ergodic problem for Hamilton-Jacobi
   equations with Neumann-type boundary conditions},
   journal={Proc. Roy. Soc. Edinburgh Sect. A},
   volume={146},
   date={2016},
   number={2},
   pages={225--242},
   issn={0308-2105},
   review={\MR{3475295}},
   doi={10.1017/S0308210515000517},
}
\bib{BaCa}{book}{
   author={Bardi, Martino},
   author={Capuzzo-Dolcetta, Italo},
   title={Optimal control and viscosity solutions of Hamilton-Jacobi-Bellman
   equations},
   series={Systems \& Control: Foundations \& Applications},
   note={With appendices by Maurizio Falcone and Pierpaolo Soravia},
   publisher={Birkh\"{a}user Boston, Inc., Boston, MA},
   date={1997},
   pages={xviii+570},
   isbn={0-8176-3640-4},
   review={\MR{1484411}},
   doi={10.1007/978-0-8176-4755-1},
}
\bib{Bar}{article}{
   author={Barles, Guy},
   title={Discontinuous viscosity solutions of first-order Hamilton-Jacobi
   equations: a guided visit},
   journal={Nonlinear Anal.},
   volume={20},
   date={1993},
   number={9},
   pages={1123--1134},
   issn={0362-546X},
   review={\MR{1216503}},
   doi={10.1016/0362-546X(93)90098-D},
}
\bib{BarB}{book}{
   author={Barles, Guy},
   title={Solutions de viscosit\'{e} des \'{e}quations de Hamilton-Jacobi},
   language={French, with French summary},
   series={Math\'{e}matiques \& Applications (Berlin) [Mathematics \&
   Applications]},
   volume={17},
   publisher={Springer-Verlag, Paris},
   date={1994},
   pages={x+194},
   isbn={3-540-58422-6},
   review={\MR{1613876}},
}
\bib{CGT}{article}{
   author={Cagnetti, Filippo},
   author={Gomes, Diogo},
   author={Tran, Hung Vinh},
   title={Adjoint methods for obstacle problems and weakly coupled systems
   of PDE},
   journal={ESAIM Control Optim. Calc. Var.},
   volume={19},
   date={2013},
   number={3},
   pages={754--779},
   issn={1292-8119},
   review={\MR{3092361}},
   doi={10.1051/cocv/2012032},
}

\bib{CLLN}{article}{
   author={Camilli, Fabio},
   author={Ley, Olivier},
   author={Loreti, Paola},
   author={Nguyen, Vinh Duc},
   title={Large time behavior of weakly coupled systems of first-order
   Hamilton-Jacobi equations},
   journal={NoDEA Nonlinear Differential Equations Appl.},
   volume={19},
   date={2012},
   number={6},
   pages={719--749},
   issn={1021-9722},
   review={\MR{2996426}},
   doi={10.1007/s00030-011-0149-7},
}

\bib{CCIZ}{article}{
   author={Chen, Qinbo},
   author={Cheng, Wei},
   author={Ishii, Hitoshi},
   author={Zhao, Kai},
   title={Vanishing contact structure problem and convergence of the viscosity solutions},
   journal={arXiv:1808.06046},
   volume={},
   date={2018},
   number={},
   pages={},
  issn={},
  review={},
   doi={},
}
\bib{CIL}{article}{
   author={Crandall, Michael G.},
   author={Ishii, Hitoshi},
   author={Lions, Pierre-Louis},
   title={User's guide to viscosity solutions of second order partial
   differential equations},
   journal={Bull. Amer. Math. Soc. (N.S.)},
   volume={27},
   date={1992},
   number={1},
   pages={1--67},
  issn={0273-0979},
  review={\MR{1118699}},
   doi={10.1090/S0273-0979-1992-00266-5},
}
\bib{CL}{article}{
   author={Crandall, Michael G.},
   author={Lions, Pierre-Louis},
   title={Viscosity solutions of Hamilton-Jacobi equations},
   journal={Trans. Amer. Math. Soc.},
   volume={277},
   date={1983},
   number={1},
   pages={1--42},
   issn={0002-9947},
   review={\MR{690039}},
   doi={10.2307/1999343},
}

\bib{DFIZ}{article}{
   author={Davini, Andrea},
   author={Fathi, Albert},
   author={Iturriaga, Renato},
   author={Zavidovique, Maxime},
   title={Convergence of the solutions of the discounted Hamilton-Jacobi
   equation: convergence of the discounted solutions},
   journal={Invent. Math.},
   volume={206},
   date={2016},
   number={1},
   pages={29--55},
   issn={0020-9910},
   review={\MR{3556524}},
   doi={10.1007/s00222-016-0648-6},
}

\bib{DZ1}{article}{
   author={Davini, Andrea},
   author={Zavidovique, Maxime},
   title={Aubry sets for weakly coupled systems of Hamilton-Jacobi
   equations},
   journal={SIAM J. Math. Anal.},
   volume={46},
   date={2014},
   number={5},
   pages={3361--3389},
   issn={0036-1410},
   review={\MR{3265180}},
   doi={10.1137/120899960},
}

\bib{DZ2}{article}{
   author={Davini, Andrea},
   author={Zavidovique, Maxime},
   title={Convergence of the solutions of discounted 
 Hamilton-Jacobi systems},
   journal={Adv. Calc. Var.},
   volume={},
   date={15 Feb 2019},
   number={},
   pages={},
   issn={},
   review={},
   doi={10.1515/acv-2018-0037},
}


\bib{EL}{article}{
   author={Engler, Hans},
   author={Lenhart, Suzanne M.},
   title={Viscosity solutions for weakly coupled systems of Hamilton-Jacobi
   equations},
   journal={Proc. London Math. Soc. (3)},
   volume={63},
   date={1991},
   number={1},
   pages={212--240},
   issn={0024-6115},
   review={\MR{1105722}},
   doi={10.1112/plms/s3-63.1.212},
}
\bib{Ev1}{article}{
   author={Evans, Lawrence C.},
   title={A survey of partial differential equations methods in weak KAM
   theory},
   journal={Comm. Pure Appl. Math.},
   volume={57},
   date={2004},
   number={4},
   pages={445--480},
   issn={0010-3640},
   review={\MR{2026176}},
   doi={10.1002/cpa.20009},
}
\bib{Ev2}{article}{
   author={Evans, Lawrence C.},
   title={Adjoint and compensated compactness methods for Hamilton-Jacobi
   PDE},
   journal={Arch. Ration. Mech. Anal.},
   volume={197},
   date={2010},
   number={3},
   pages={1053--1088},
   issn={0003-9527},
   review={\MR{2679366}},
   doi={10.1007/s00205-010-0307-9},
}
		
\bib{Fa1}{article}{
   author={Fathi, Albert},
   title={Th\'{e}or\`eme KAM faible et th\'{e}orie de Mather sur les syst\`emes
   lagrangiens},
   language={French, with English and French summaries},
   journal={C. R. Acad. Sci. Paris S\'{e}r. I Math.},
   volume={324},
   date={1997},
   number={9},
   pages={1043--1046},
   issn={0764-4442},
   review={\MR{1451248}},
   doi={10.1016/S0764-4442(97)87883-4},
}
	
\bib{Fa2}{book}{
   author={Fathi, Albert},
   title={Weak KAM theorem in Lagrangian dynamics, Preliminary version 10},
   language={English},
   series={},
   volume={},
   publisher={},
   date={2008},
   pages={xiv+273},
   isbn={},
   review={},
}


\bib{Go}{article}{
   author={Gomes, Diogo Aguiar},
   title={Duality principles for fully nonlinear elliptic equations},
   conference={
      title={Trends in partial differential equations of mathematical
      physics},
   },
 book={
      series={Progr. Nonlinear Differential Equations Appl.},
      volume={61},
      publisher={Birkh\"{a}user, Basel},
   },
   date={2005},
   pages={125--136},
   review={\MR{2129614}},
   doi={10.1007/3-7643-7317-2\_10},  
}
\bib{GMT}{article}{
   author={Gomes, Diogo A.},
   author={Mitake, Hiroyoshi},
   author={Tran, Hung V.},
   title={The selection problem for discounted Hamilton-Jacobi equations:
   some non-convex cases},
   journal={J. Math. Soc. Japan},
   volume={70},
   date={2018},
   number={1},
   pages={345--364},
   issn={0025-5645},
   review={\MR{3750279}},
   doi={10.2969/jmsj/07017534},
}

 
\bib{IJ}{article}{
   author={Ishii, Hitoshi},
   author={Jin, Liang},
   title={The vanishing discount problem for 
monotone systems of Hamilton-Jacobi equations.  Part 2: Nonlinear coupling},
   journal={Calc. Var. Partial Differential Equations},
   volume={},
   date={},
   number={},
   pages={},
   issn={},
  note={(to appear)},
   review={},
   
   doi={10.1007/s00526-020-01768-8},
}


\bib{IPerron}{article}{
   author={Ishii, Hitoshi},
   title={Perron's method for Hamilton-Jacobi equations},
   journal={Duke Math. J.},
   volume={55},
   date={1987},
   number={2},
   pages={369--384},
   issn={0012-7094},
   review={\MR{894587}},
   doi={10.1215/S0012-7094-87-05521-9},
}

\bib{IK}{article}{
   author={Ishii, Hitoshi},
   author={Koike, Shigeaki},
   title={Viscosity solutions for monotone systems of second-order elliptic
   PDEs},
   journal={Comm. Partial Differential Equations},
   volume={16},
   date={1991},
   number={6-7},
   pages={1095--1128},
   issn={0360-5302},
   review={\MR{1116855}},
   doi={10.1080/03605309108820791},
}

\bib{IMT1}{article}{
   author={Ishii, Hitoshi},
   author={Mitake, Hiroyoshi},
   author={Tran, Hung V.},
   title={The vanishing discount problem and viscosity Mather measures. Part
   1: The problem on a torus},
   language={English, with English and French summaries},
   journal={J. Math. Pures Appl. (9)},
   volume={108},
   date={2017},
   number={2},
   pages={125--149},
   issn={0021-7824},
   review={\MR{3670619}},
   doi={10.1016/j.matpur.2016.10.013},
}
\bib{IMT2}{article}{
   author={Ishii, Hitoshi},
   author={Mitake, Hiroyoshi},
   author={Tran, Hung V.},
   title={The vanishing discount problem and viscosity Mather measures. Part
   2: Boundary value problems},
   language={English, with English and French summaries},
   journal={J. Math. Pures Appl. (9)},
   volume={108},
   date={2017},
   number={3},
   pages={261--305},
   issn={0021-7824},
   review={\MR{3682741}},
   doi={10.1016/j.matpur.2016.11.002},
}
\bib{IS}{article}{
   author={Ishii, Hitoshi},
   author={Siconolfi, Antonio},
   title={The vanishing discount problem for Hamilton--Jacobi equations in
   the Euclidean space},
   journal={Comm. Partial Differential Equations},
   volume={45},
   date={2020},
   number={6},
   pages={525--560},
   issn={0360-5302},
   review={\MR{4106998}},
   doi={10.1080/03605302.2019.1710845},
}

\bib{MSTY}{article}{
   author={Mitake, H.},
   author={Siconolfi, A.},
   author={Tran, H. V.},
   author={Yamada, N.},
   title={A Lagrangian approach to weakly coupled Hamilton-Jacobi systems},
   journal={SIAM J. Math. Anal.},
   volume={48},
   date={2016},
   number={2},
   pages={821--846},
   issn={0036-1410},
   review={\MR{3466199}},
   doi={10.1137/15M1010841},
}

\bib{MT}{article}{
   author={Mitake, Hiroyoshi},
   author={Tran, Hung V.},
   title={Selection problems for a discount degenerate viscous
   Hamilton-Jacobi equation},
   journal={Adv. Math.},
   volume={306},
   date={2017},
   pages={684--703},
   issn={0001-8708},
   review={\MR{3581314}},
   doi={10.1016/j.aim.2016.10.032},
}
\bib{MT1}{article}{
   author={Mitake, Hiroyoshi},
   author={Tran, Hung V.},
   title={Remarks on the large time behavior of viscosity solutions of
   quasi-monotone weakly coupled systems of Hamilton-Jacobi equations},
   journal={Asymptot. Anal.},
   volume={77},
   date={2012},
   number={1-2},
   pages={43--70},
   issn={0921-7134},
   review={\MR{2952714}},
}


\bib{MT2}{article}{
   author={Mitake, H.},
   author={Tran, H. V.},
   title={A dynamical approach to the large-time behavior of solutions to
   weakly coupled systems of Hamilton-Jacobi equations},
   language={English, with English and French summaries},
   journal={J. Math. Pures Appl. (9)},
   volume={101},
   date={2014},
   number={1},
   pages={76--93},
   issn={0021-7824},
   review={\MR{3133425}},
   doi={10.1016/j.matpur.2013.05.004},
}

\bib{MT3}{article}{
   author={Mitake, Hiroyoshi},
   author={Tran, Hung V.},
   title={Homogenization of weakly coupled systems of Hamilton-Jacobi
   equations with fast switching rates},
   journal={Arch. Ration. Mech. Anal.},
   volume={211},
   date={2014},
   number={3},
   pages={733--769},
   issn={0003-9527},
   review={\MR{3158806}},
   doi={10.1007/s00205-013-0685-x},
}
		
\bib{PLL}{book}{
   author={Lions, Pierre-Louis},
   title={Generalized solutions of Hamilton-Jacobi equations},
   series={Research Notes in Mathematics},
   volume={69},
   publisher={Pitman (Advanced Publishing Program), Boston, Mass.-London},
   date={1982},
   pages={iv+317},
   isbn={0-273-08556-5},
   review={\MR{667669}},
}

\bib{LPV}{article}{
  author={Lions,  P.-L.}, 
  author={Papanicolaou, G.},
  author={Varadhan, S.}, 
  title={Homogenization of Hamilton-Jacobi equations}, 
  journal={unpublished work},
  date={1987},
}

\bib{Si}{article}{
   author={Sion, Maurice},
   title={On general minimax theorems},
   journal={Pacific J. Math.},
   volume={8},
   date={1958},
   pages={171--176},
   issn={0030-8730},
   review={\MR{0097026}},
}

\bib{Ter}{article}{
   author={Terai, Kengo},
   title={Uniqueness structure of weakly coupled systems of
ergodic problems of Hamilton-Jacobi equations},
   journal={ArXiv:1901.05314v1},
   volume={},
   date={},
   pages={},
   issn={},
   review={},
   doi={},
}

\bib{Te}{article}{
   author={Terkelsen, Frode},
   title={Some minimax theorems},
   journal={Math. Scand.},
   volume={31},
   date={1972},
   pages={405--413 (1973)},
   issn={0025-5521},
   review={\MR{0325880}},
   doi={10.7146/math.scand.a-11441},
}


\bib{Var}{book}{
   author={Varga, Richard S.},
   title={Matrix iterative analysis},
   series={Springer Series in Computational Mathematics},
   volume={27},
   edition={Second revised and expanded edition},
   publisher={Springer-Verlag, Berlin},
   date={2000},
   pages={x+358},
   isbn={3-540-66321-5},
   review={\MR{1753713}},
   doi={10.1007/978-3-642-05156-2},
}

\end{biblist}
\end{bibdiv}
\bye